\documentclass[11pt,fleqn]{article}
\usepackage{amsfonts}
\newcommand{\ol}{\setlength{\itemsep}{0pt.}\begin{enumerate}}
\newcommand{\eol}{\end{enumerate}\setlength{\itemsep}{-\parsep}}
\newcommand{\ignore}[1]{}
\title{Low-degree tests at large distances}
\author{Alex Samorodnitsky\thanks{Institute of Computer Science,
    Hebrew University. 
{\tt salex@huji.ac.il} This paper
is based upon work supported by the Israel Science Foundation under
grant 039-716 and by the German-Israeli Foundation under grant I-2052.}}

\begin{document}

\begin{titlepage}

\maketitle

\sloppy

\thispagestyle{empty}

\begin{abstract}
We define tests of boolean functions which
distinguish between linear (or quadratic) polynomials, and functions
which are very far, in an appropriate sense, from these 
polynomials. The tests have optimal or nearly optimal trade-offs
between soundness and the number of queries.

In particular, we show that functions with small Gowers uniformity
norms behave ``randomly'' with respect to hypergraph linearity tests.

A central step in our analysis of quadraticity tests is the proof of an
inverse theorem for the third Gowers uniformity norm of boolean functions.

The last result has also a coding theory application. It is
possible to estimate efficiently the distance from the second-order
Reed-Muller code on inputs lying far beyond its list-decoding radius.

\end{abstract}

\end{titlepage}


\newtheorem{THEOREM}{Theorem}[section]
\newenvironment{theorem}{\begin{THEOREM} \hspace{-.85em} {\bf :}
}%
                        {\end{THEOREM}}
\newtheorem{LEMMA}[THEOREM]{Lemma}
\newenvironment{lemma}{\begin{LEMMA} \hspace{-.85em} {\bf :} }%
                      {\end{LEMMA}}
\newtheorem{COROLLARY}[THEOREM]{Corollary}
\newenvironment{corollary}{\begin{COROLLARY} \hspace{-.85em} {\bf
:} }%
                          {\end{COROLLARY}}
\newtheorem{PROPOSITION}[THEOREM]{Proposition}
\newenvironment{proposition}{\begin{PROPOSITION} \hspace{-.85em}
{\bf :} }%
                            {\end{PROPOSITION}}
\newtheorem{DEFINITION}[THEOREM]{Definition}
\newenvironment{definition}{\begin{DEFINITION} \hspace{-.85em} {\bf
:} \rm}%
                            {\end{DEFINITION}}
\newtheorem{EXAMPLE}[THEOREM]{Example}
\newenvironment{example}{\begin{EXAMPLE} \hspace{-.85em} {\bf :}
\rm}%
                            {\end{EXAMPLE}}
\newtheorem{CONJECTURE}[THEOREM]{Conjecture}
\newenvironment{conjecture}{\begin{CONJECTURE} \hspace{-.85em}
{\bf :} \rm}%
                            {\end{CONJECTURE}}
\newtheorem{MAINCONJECTURE}[THEOREM]{Main Conjecture}
\newenvironment{mainconjecture}{\begin{MAINCONJECTURE} \hspace{-.85em}
{\bf :} \rm}%
                            {\end{MAINCONJECTURE}}
\newtheorem{PROBLEM}[THEOREM]{Problem}
\newenvironment{problem}{\begin{PROBLEM} \hspace{-.85em} {\bf :}
\rm}%
                            {\end{PROBLEM}}
\newtheorem{QUESTION}[THEOREM]{Question}
\newenvironment{question}{\begin{QUESTION} \hspace{-.85em} {\bf :}
\rm}%
                            {\end{QUESTION}}
\newtheorem{REMARK}[THEOREM]{Remark}
\newenvironment{remark}{\begin{REMARK} \hspace{-.85em} {\bf :}
\rm}%
                            {\end{REMARK}}

\newcommand{\thm}{\begin{theorem}}
\newcommand{\lem}{\begin{lemma}}
\newcommand{\pro}{\begin{proposition}}
\newcommand{\dfn}{\begin{definition}}
\newcommand{\rem}{\begin{remark}}
\newcommand{\xam}{\begin{example}}
\newcommand{\cnj}{\begin{conjecture}}
\newcommand{\mcnj}{\begin{mainconjecture}}
\newcommand{\prb}{\begin{problem}}
\newcommand{\que}{\begin{question}}
\newcommand{\cor}{\begin{corollary}}
\newcommand{\prf}{\noindent{\bf Proof:} }
\newcommand{\ethm}{\end{theorem}}
\newcommand{\elem}{\end{lemma}}
\newcommand{\epro}{\end{proposition}}
\newcommand{\edfn}{\bbox\end{definition}}
\newcommand{\erem}{\bbox\end{remark}}
\newcommand{\exam}{\bbox\end{example}}
\newcommand{\ecnj}{\bbox\end{conjecture}}
\newcommand{\emcnj}{\bbox\end{mainconjecture}}
\newcommand{\eprb}{\bbox\end{problem}}
\newcommand{\eque}{\bbox\end{question}}
\newcommand{\ecor}{\end{corollary}}
\newcommand{\eprf}{\bbox}
\newcommand{\beqn}{\begin{equation}}
\newcommand{\eeqn}{\end{equation}}
\newcommand{\blist}{\begin{itemize}}
\newcommand{\elist}{\end{itemize}}
\newcommand{\wbox}{\mbox{$\sqcap$\llap{$\sqcup$}}}
\newcommand{\bbox}{\vrule height7pt width4pt depth1pt}
\newcommand{\qed}{\bbox}
\def\sup{^}

\def\B{\{0,1\}} 
\def\H{\{-1,1\}}

\def\S{S(n,w)}

\def\n{\lfloor \frac n2 \rfloor}

\def\Tp{Tchebyshef polynomial}
\def\Tps{TchebysDeto be the maximafine $A(n,d)$ l size of a code with
  distance $d$hef polynomials} 
\newcommand{\rarrow}{\rightarrow}

\newcommand{\larrow}{\leftarrow}

\overfullrule=0pt
\def\setof#1{\lbrace #1 \rbrace}

\def \E{\mathbb E}
\def \R{\mathbb R}
\def \Z{\mathbb Z}
\def \F{\mathbb F}

\def\<{\left<}
\def\>{\right>}
\def \({\left(}
\def \){\right)}
\def \[{\left[}
\def \]{\right]}
\def \e{\epsilon}
\def \d{\delta}
\def \l{\lambda}
%
\section{Introduction}
This paper returns to the general question of the relation between number of
queries and the probability of error in low-degree tests. 

The specific questions we deal with originate within a wider
framework of Probabilistically Checkable Proofs (PCPs). The PCP
theorem \cite{ALMSS92,AS92} states that it is possible to encode
certificates of 
satisfiability for SAT instances in such a way that a probabilistic
verifier using logarithmic number of random bits can check the
validity of the certificate with high probability of success, after
looking only at a constant number of bits in the encoding. We 
consider here only PCPs with almost perfect {\it completeness}, which
means that valid certificates are nearly always
\footnote{See, e.g., \cite{Trevisan98}  for a precise definition}
accepted. Given this, 
and fixing the number of queries $q$, we are interested in the best
possible {\it soundness} of the PCP, namely the probability $s$ of accepting
an encoding of a false proof. 

It is easy to see that, unless $P = NP$, the lower bound $s \ge 1/2^q$
must hold. Stronger lower bounds were given in \cite{H,Trevisan96}. 
The best known lower bound \cite{CMM} is
$s \ge \Omega\(\frac{q}{2^q}\)$. From the other direction, the PCP
theorem shows that we can achieve $s \le \frac{1}{2^{O(q)}}$, and it
was shown in \cite{EH}, following \cite{ST00}, that $s \le
\frac{2^{\sqrt{2q}}}{2^q}$. In 
  \cite{ST06}, assuming the Unique Games Conjecture \cite{K}, the
  upper bound was improved to $s \le q/2^{q-1}$, which is of
  course (conditionally) best possible, up to constants.   

Let us say a few words on the structure of a PCP protocol. In the
common paradigm \cite{AS92} the verifier of a PCP is split into two
entities, the {\it inner} and the 
{\it outer} verifiers. Roughly speaking, the outer verifier chooses the
(randomized) portion of the proof to be checked by the inner verifier. The
inner verifier views the binary string it is given as a boolean
function, and looks for a certain combinatorial pattern. If the
pattern is not there the proof is rejected. If the 
inner verifier finds the appropriate property with non-negligible
probability over its inputs, the outer verifier can then use this
information to validate the PCP statement.

Due to the gap structure inherent in the PCP construction, the decision
of the inner verifier is usually dichotomic. This is to say it must accept
if the property is satisfied and reject only if the function is very far,
in the appropriate sense, from having the property. 

In this framework, an often considered property of a boolean function
is that of being represented by a low-degree polynomial over a finite
field. Here we deal only with the field of two elements and this
representation is particularly simple: 
\dfn
A boolean function $f:\{0,1\}^n \rarrow \{-1,1\}$ has a degree-$d$
representation if $f(x) = (-1)^{P(x)}$, where $P(x) = P\(x_1...x_n\)$
is an $n$-variate polynomial of degree $d$ over $\F_2$.  
\edfn
In our version of the Low-Degree testing problem we are given an
oracle access to a boolean function $f:\{0,1\}^n \rarrow \{-1,1\}$ and
we want to determine whether 
\begin{enumerate}
\item
The function $f$ can be represented by a degree-$d$ polynomial
\item
It is $\frac12 - \e$ far from any function with such representation. 
\end{enumerate}
The distance between two functions is a fraction of points in which
they disagree. 

Low-degree tests we consider have perfect completeness, namely
in case (1) they always accept. We now define the soundness of a test.
\dfn
\label{def:soundness}
A low-degree test has soundness $s$ if for any function
$f$ that is $\frac12 - \epsilon$ far from degree-$d$ polynomials,
the test accepts $f$ with probability at most $s + \phi(\epsilon)$
where $\phi(x) \rarrow_{x\rarrow 0} 0$. 
\edfn
Designing a low-degree test with a good trade-off between the number
of queries and the soundness is a step towards a PCP construction. There
are several ways in which such a result needs to be augmented to lead
to a full PCP construction. We refer to the discussion in
\cite{ST06}. It seems, however, that in most cases in which this
extension process succeeded, the obtained PCP inherited the relevant
parameters (number of queries, soundness) of the low-degree test. 

Degree-$1$ (linear) tests with asymptotically optimal asymptotic
trade-off between the 
number of queries and the soundness where given in \cite{ST00}. In the
same paper these tests were extended to PCP constructions with
similar parameters.  

A natural way to improve the PCP parameters further is to
consider additional combinatorial tests.   

In this paper we study degree-$2$ tests and relaxed
versions of the degree-$1$ test. Here is a brief overview of our main
results.
\begin{itemize}
\item
We define and analyze a degree-$1$ test with relaxed rejection criteria whose
trade-off between the number of queries and the soundness is
asymptotically optimal and is much
better than that achievable by the standard linearity tests. A
different (and easier) analysis of this test was 
given in \cite{ST06}. In that paper we were also able to extend the
test to a conditional PCP construction (assuming the Unique Games
Conjecture \cite{K}) with an optimal  number of queries vs. soundness
trade-off.  
\item
We define and analyze a degree-$2$ test with a very
good trade-off between the number of queries $q$ and the soundness
$s$. (We conjecture this trade-off to be asymptotically optimal.) A
technical ingredient of this 
result has a natural interpretation in the framework of error-correcting
codes. We give a tight analysis of the acceptance
probability of a natural local test of \cite{AKKLR} for the second-order
Reed-Muller code at distances near the covering radius of this
code. As a consequence, it turns out to be possible to estimate
efficiently the distance from this code on inputs lying far beyond its
list-decoding radius. 
\end{itemize}
Our analysis of these tests is based on several technical assertions which
could be of independent interest, and which we describe next.
\begin{itemize}
\item
We give a tight analysis of the Abelian Homomorphism testing problem
for some families of groups, including powers of $\Z_p$. The central
technical claim, which we state here for the special case of $p=2$, is
that if a function 
$\phi:~\Z^n_2 \rarrow \Z^n_2$ satisfies $Pr\(\phi(x+y) = \phi(x) +
\phi(y)\)$ with probability bounded away from zero, then there is a
matrix $D\in M_{n,n}\(\Z_2\)$ such that a linear transformation
$\psi:x \mapsto Dx$ coincides with $\phi$ on a non-negligible fraction
of the inputs.
\item
We introduce and study the notion of a {\it generalized average} of a
function $f$ over $\B^n$. A generalized average is a non-linear
functional on the space of real (or complex) valued functions on the
boolean cube. It is associated with a binary
matrix $M$ and it measures the average over a certain family of
subsets of $\B^n$, defined by $M$, of products of $f$ over each
subset. Generalized averages arize naturally in the analysis of 
low-degree tests. 
An important special case is when this family consists of
all the affine subsets of $\B^n$ of a fixed dimension $d$. The
generalized average in this case turns out to measure (a power of) a
norm of the function $f$. These norms are the {\it Gowers uniformity
  norms} \cite{Gowers01} and they measure, in a certain sense, 
a proximity of the function to a polynomial of degree $d$. 
\begin{itemize}
\item
We show that a function with a large third uniformity norm is somewhat
close to an $n$-variate quadratic polynomial over $\F_2$. Similar
results for finite Abelian groups of cardinality indivisible by $6$ have been
independently proved in \cite{GT05}. 
\item
We show that functions on which the {\it hypergraph linearity tests}
defined in \cite{ST00} fail with non-negligible probability have large
uniformity norms. 
\item
We observe that functions with small uniformity norms are {\it
  pseudorandom} in the sense of \cite{Gowers01}, and briefly discuss
  pseudorandom properties of such functions in our context.
\end{itemize}
\end{itemize}

In the next sections we give a more detailed description of the
background and of the results in this paper. The proofs are given in
the Appendices.  

\noindent {\it Organization}

We describe relaxed linearity tests in
Section~\ref{sec:deg1}. Degree-$2$ tests and properties
of the Reed-Muller code of order $2$ are described in
Section~\ref{sec:deg2}. Abelian homomorphism testing is discussed in
Section~\ref{sec:hom}. Section~\ref{sec:technical} gives more details
on the technical tools used, in particular their connections with
recent work in additive number theory. A notion of pseudorandomness of
boolean functions which comes from additive number theory is
introduced and briefly discussed.

\section{Degree-$1$ Tests}
\label{sec:deg1}
This is the simplest and the most useful case in practice. A boolean
function $f$ has a degree-$1$ representation if $f(x) = (-1)^{\<a,x\>
  + b}$, where $a$ is a fixed vector in $\B^n$, and $b$ is a fixed
constant in $\B$. 
Hence in this case the tester has to decide whether the function is
linear
\footnote{or rather {\it affine}. In practice the function is usually
  tested for linearity ($b=0$). The two testing problems are
  essentially equivalent, and we occasionally 
  will, with some abuse of meaning, refer to both as {\it linearity
  testing problems}.} or is far from every linear function.

A simple linearity test with three queries was defined in \cite{BLR}. 
\footnote{We observe that to transform this test to an affinity
  (degree-$1$) test, it 
suffices to replace $1$ with $f(0)$ in the definition of the test.} 
\smallskip

\noindent\fbox{\begin{minipage}{3.0in}
Choose uniformly at random $x,y \in \B^n$\\
If $f(x)f(y)f(x+y) = 1$\\ 
then accept\\
else reject
\end{minipage}}

\medskip
It is shown in \cite{johan+} that if this test accepts $f$ with
probability $\frac12 + \delta$ then $f$ is $\frac12 - 2\delta$
close to a linear function. Therefore, according to our definition,
this test has soundness $s = \frac12$.

Independent repetition of $q/3$ basic tests leads to a test with $q$
queries and soundness $s = \frac{1}{2^{2q/3}}$. 
To improve the trade-off between $q$ and $s$, more complex tests have
to be considered. It turns out that it is possible to associate such a
test with any given graph. Fix a graph $G=(V,E)$ on $t$ vertices, The
following test is a {\it dependent} combination of the basic tests of
\cite{BLR}. 
\smallskip

\noindent\fbox{\begin{minipage}{2.5in}
 Choose uniformly at random $x_1,\ldots,x_{t} \in \B^n$\\
 If $\prod_{i\in e} f(x_i) \cdot f(\sum_{i\in e} x_i)=f(0)$ for
 all $e \in E$ \\ 
 then accept\\
 else reject
\end{minipage}}

\medskip
These {\it graph tests} were defined in \cite{Trevisan98}. A graph test
associated with a graph $G = (V,E)$ runs $|E|$ correlated copies of
the basic linearity test. 
In \cite{ST00} it was shown that for functions which are far from
degree-$1$ polynomials (this is to say, have small Fourier
coefficients), these copies of the basic test behave 
essentially independently. More precisely, the soundness of this test
is $s = 1/2^{|E|}$. Of course, the total number of queries is $q = |V|
+ |E|$. In particular, choosing $G$ to be the complete graph on $t$
vertices, we obtain an affinity test with $q = {t\choose 2} + t$ and
$s = \frac{1}{2^{{t \choose 2}}}$. This means that $s \approx
  \frac{2^{\sqrt{2q}}}{2^q}$. 

A natural generalization of graph tests to {\it hypergraph tests} was
given in \cite{ST00}. Let $H = (V,E)$ be a hypergraph on $t$
vertices and consider the following test:
\smallskip

\noindent\fbox{\begin{minipage}{2.5in}
 Choose uniformly at random $x_1,\ldots,x_{t} \in \B^n$\\
 If $\prod_{i\in e} f(x_i) \cdot f(\sum_{i\in e} x_i)=f^{|e|+1}(0)$ for
 all $e \in E$ \\ 
 then accept\\
 else reject
\end{minipage}}

\medskip
A hypergraph test runs $|E|$ copies of the basic linearity test, where
$|E|$ is now the number of hyper-edges. Unfortunately, it is not true
that, for functions far from degree-$1$ 
polynomials, these copies behave independently. Consider a function
$f(x) = (-1)^{x_1 x_2 + ... + x_{n-1} x_n}$. This function is
maximally far from all degree-$1$ polynomials (it is a {\it bent
  function}), but any hypergraph test with $q = |V| + |E|$ queries
accepts this function with probability at least
$\frac{2^{\Omega\(\sqrt q\)}}{2^q}$ \cite{ST00}. More generally, we show in
\cite{ST06} that this is true for any non-adaptive linearity test that
always accepts linear functions.

\noindent {\bf Our results}

The starting point of this work was the realization that the function $f$
we have described is a quadratic polynomial and that it is accepted by
a hypergraph test with non-negligible probability, because, roughly
speaking, the basic ingredient of this test takes a 
discrete derivative of the tested function and compares it to
zero. The order of the derivative is essentially given by the
cardinality of the hyperedges. We will say more about this in
Section~\ref{app:A} and in the full version of the paper. The
natural question then is whether quadratic polynomials, and, more
generally, low-degree polynomials, are the only obstructions to better
performance by hypergraph linearity tests. 

We give a partial (affirmative) answer to this question for general
hypergraphs. We are able to answer this question completely for hypergraphs of
maximal edge-size $3$ and for quadratic polynomials. The answer is again
positive. We conjecture the answer to be positive for general hypergraphs and 
low-degree polynomials. 

We prove two claims. These are the main technical results of this
paper.

The first claim is valid for any hypergraph. First, we define Gowers
uniformity norms.  
\dfn
\label{def:gowers}
Let $f:~\B^n \rarrow \R$ be a function, and $d \ge 1$ be an
integer. The $d$-th Gowers uniformity norm (for the group $\Z^n_2$)
is given by
$$
\|f\|_{U_d} = \left[ \E_{x,y_1,...,y_d} \prod_{S \subseteq [d]} f\(x +
\sum_{i\in S} y_i\) \right]^{1/2^d}
$$
Here $x,y_1,...,y_d$ are chosen uniformly and independently at random 
from $\B^n$.
\edfn
\thm
\label{thm:gowers-norms}
Let $H = (V,E)$ be a hypergraph with maximal edge-size $d$. Then the
probability that the linearity test associated with $H$ accepts a
boolean function $f$ is bounded by
$$
\frac{1}{2^{|E|}} + \|f\|_{U_d}
$$
\ethm
Another (and easier) proof of this theorem and its generalization to
several functions is given in \cite{ST06}.

The second claim is that a boolean function with a large third
unformity norm is somewhat close to a quadratic polynomial. 
\thm
\label{thm:inverse-quadratic}
Let $f:~\B^n \rarrow \H$ be a function such that $\|f\|_{U_3} \ge
\e$. Then there exists a quadratic polynomial $g$ such that the
distance between $f$ and $g$ is at most $\frac12 - \e'$. Here one can choose
$\e' \ge \Omega\(\mbox{exp}\left\{-\(\frac{1}{\e^C}\)\right\}\)$ for 
an absolute constant $C$.
\ethm

\ignore{Theorem~\ref{thm:gowers-norms} leads directly to a definition of a
relaxed linearity test with an essentially optimal trade-off between
the number of queries and the soundness. 
}
Consider the following relaxed degree-$1$ testing problem. Given an
oracle access to a boolean function $f:\{0,1\}^n \rarrow \{-1,1\}$ and
an integer $d \ge 2$ we want to determine whether 
\begin{enumerate}
\item
The function $f$ can be represented by a degree-$1$ polynomial.
\item
$\|f\|_{U_d} \le \e$.
\end{enumerate}
Once again we want tests with perfect completeness. The soundness of
the test is defined as in Definition~\ref{def:soundness}.
\rem
Let us point out, that this test is indeed a relaxation of the standard'
degree-$1$ test. It is known \cite{Gowers01} that uniformity norms
$\|f\|_{U_d}$ of $f$ 
are monotone increasing in $d$. It is easy to see that the second
uniformity norm is the same 
as the $l_4$ norm of the Fourier transform of $f$: $\|f\|_{U_2} =
\(\sum_{\alpha \in \B^n} \hat{f}^4(\alpha)\)^{\frac14} \ge
\max_{\alpha \in \B^n} |\hat{f}(\alpha)|$. This means that the
functions the test has to reject are at least $\frac12 - \frac{\e}{2}$ far from
degree-$1$ polynomials.
\erem
It is a direct consequence of Theorem~\ref{thm:gowers-norms} that
hypergraph tests solve the relaxed testing problem with the ``right''
soundness.
\thm
\label{thm:hypergraph-test}
Let $d\ge 2$ and let $H = (V,E)$ be a hypergraph with maximal
edge-size $d$. Then the hypergraph linearity test associated with $H$
solves the relaxed 
degree-$1$ testing problem with perfect completeness and soundness
$1/2^{|E|}$. 
\ethm
Choosing $H$ to be a complete $d$-uniform hypergraph on $t \approx
q^{1/d}$ vertices
leads to a test with $q$ queries and soundness $s \le
\frac{2^{\Omega\(q^{1/d}\)}}{2^q}$. 
This trade-off is shown to be asymptotically optimal in \cite{ST06}.

It remains to observe that Theorem~\ref{thm:hypergraph-test} together
with Theorem~\ref{thm:inverse-quadratic} imply that the complete 
$3$-uniform hypergraph test distinguishes between linear functions and
functions which are far from quadratic polynomials with optimal soundness of 
$s \le \frac{2^{\Omega\(q^{1/3}\)}}{2^q}$.
\section{Second-Order Reed-Muller Codes}
\label{sec:deg2}
A binary error-correcting code \cite{mcw-sloane} of length $N$ and
normalized distance 
$\d$ is a subset of $\B^N$ in which any two distinct elements disagree
on at least $\d$-fraction of the domain (the coordinates). This allows
for error-correction: a corrupted codeword (element of the code) with
less than $\d/2$-fraction of the errors can, in principle, be 
recovered by going to
the unique nearest element of the code. We call $\d/2$ the
unique-decoding radius of the code. 

Finding the nearest codeword
can be computationally hard. Here we are interested in efficient
error-correction. 

An important example of a code of length $N = 2^n$ is the subset of $\B^N$
whose elements are evaluations of $n$-variate degree $d$ polynomials 
over $\F_2$. This is
the {\it Reed-Muller (RM) code} of order $d$. Efficient error-correcting
algorithms for RM codes were given in \cite{Reed}. 

One can go beyond unique decoding. It is an easy consequence of
the Johnson bound for constant-weight codes \cite{mcw-sloane} that 
there is $\l > \d/2$ such that there could be only a few (polynomially
many in $N$) 
codewords at distance $\l$ from a corrupted codeword. We call maximal
$\l$ with this property the list-decoding radius of the
code. For many codes 
there are efficient list-decoding algorithms \cite{Madhu_survey} that,
for any $\l$ smaller than list-decoding radius,
recover all the codewords within distance $\l$ from the corrupted
codeword. To the best
of our knowledge there are no such algorithms for binary RM codes of
order larger than $1$.

Another useful property of a code is local testability
\cite{Goldreich-Sudan}. A code is locally testable if 
there exists an efficient randomized algorithm (test) which, given an access
to a putative codeword $f \in \B^N$, examines a finite
number of coordinates of $f$ and decides whether
$f$ is a codeword. We want the test always to accept valid codewords,
and to minimize the probability $s$ of accepting an invalid codeword,
given the number of queries $q$. The questions we discuss in this
paper fall naturally into the framework of local testability of
Reed-Muller codes. In fact, we deal with a special case (a promise
problem) in which the putative codeword is promised either to lie in
the code or to be ($\frac12 - \e$) far from the code. We remark that, in the
general case, the probability the test accepts an invalid codeword 
will necessarily depend also on its distance from the code.

\xam
A good example for the notions we have discussed is the first order
Reed-Muller code, also known as the Hadamard code. The
distance of this code is $1/2$, and therefore its unique-decoding
radius is $1/4$. However, it is efficiently list-decodable for any
distance $\l < 1/2$ \cite{Goldreich-Levin}. 

The Hadamard code is also locally testable. In fact, the basic
linearity test of \cite{BLR} is a good $3$-query local test for this code. 
\cite{johan+} studies the dependence of the probability this test accepts
an invalid codeword on its distance from the code. For
distances close to $1/2$ the analysis it tight, 
and the probability of acceptance is shown to be upper bounded by $1$
minus the distance. 
\exam

Local testability of Reed-Muller codes of any fixed order $d$ was proved
in \cite{AKKLR}. The basic test in \cite{AKKLR} (presented here with a
small twist to adopt it to our setting) chooses independently at random
$d+2$ vectors $x,y_1,...,y_{d+1}$ in $\B^n$, and computes the product of
the tested function over the $d$-dimensional affine subspace of $\B^n$
given by $x + Span\(y_1,...,y_{d+1}\)$. If the product is $1$ the test
accepts. Otherwise it rejects. This is a natural generalization
of the linearity test of \cite{BLR}. While that test can be
interpeted as taking a random second directional derivative
and checking whether it vanishes, the test of \cite{AKKLR} amounts to 
checking whether a random 
derivative of order $d+1$ vanishes. \cite{AKKLR} studies the
dependence of the probability this test accepts 
an invalid codeword $f$ on its distance from the code. (We
observe that this probability is precisely $\frac{1 +
  \|f\|^{2^d}_{U_d}}{2}$, cf. Definition~\ref{def:gowers}). In 
particular it is shown that, for distances larger than $2^{-d}$, the
probability of acceptance is upper bounded by $1 - \Omega\(d^{-1}
2^{-d}\)$. Thus, for $d=2$, the probability of acceptance is
upper bounded by some constant smaller than $1$.

\noindent {\bf Our results}

We study the probability of error of the test of \cite{AKKLR} for the
second-order Reed-Muller code and for
distances close to $1/2$. We provide a tight analysis for this case,
showing this probability to be essentially
upper bounded by $1$ minus the
distance. Specifically, by Theorem~\ref{thm:inverse-quadratic},
if this probability is larger that $1/2 + \e$ then there is a
quadratic polynomial whose distance from the tested function is at
most $1/2 - \e'$.

Our result has a following coding interpretation. Although the 
list-decoding raduis of the second-order Reed-Muller code is $1/4$
\cite{mcw-sloane}, it
is possible to determine whether the distance of a given function $f
\in \B^N$ from the code is strictly smaller than 
the covering radius of the code, which is $1/2 - o(1)$.~\footnote{This is also,
with overwhelming probability, 
the typical distance of an element of $\B^N$ from the code.} More
precisely, we have the following proposition. 
\pro
\label{pro:dicho}
There is a positive constant $C$ such that, given a
function $f:~\B^n \rarrow \H$, and a parameter $\d > 0$, 
it is possible to determine, with probability arbitrarily close to
$1$, and 
in time linear in $\frac{1}{\d}$, which of the two following
(mutually non-exclusive) options holds:
\begin{itemize}
\item
The distance of $f$ from the quadratic polynomials is at least
$\frac12 - \Omega\(\d^{1/2}\)$.
\item
The distance of $f$ from the quadratic polynomials is at most
$\frac12 - \exp\left\{\(-\frac{1}{\d}\)^C\right\}$.
\end{itemize}
\epro
Combining theorems~\ref{thm:hypergraph-test} and
\ref{thm:inverse-quadratic} leads to our main result in this
section, an analysis of hypergraph degree-$2$ (quadraticity) tests. 

Given a $3$-uniform hypergraph
$H=([t],E)$ on $t$ vertices, the test is defined as follows.
\smallskip

\noindent\fbox{\begin{minipage}{5in}
 Choose uniformly at random $x_1,\ldots,x_{t} \in \B^n$\\
 If for
 all $e = \{i,j,k\}
 \in E$ holds\\
${\E}_{x_i,x_j,x_k}
 f(x_i)f(x_j)f(x_k)f(x_i+x_j)f(x_i+x_k)f(x_j+x_k)f(x_i+x_j+x_k) = f(0)
$
 then accept\\
 else reject
\end{minipage}}

\medskip
\thm
\label{thm:hypergraph-quadraticity-test}
Let $H = (V,E)$ be a $3$-uniform hypergraph. 
Then the hypergraph quadraticity test solves the 
degree-$2$ testing problem with perfect completeness and soundness
$1/2^{|E|}$. 
\ethm
Choosing $H$ to be a complete $3$-uniform hypergraph on $t \approx
q^{1/3}$ vertices
leads to a test with $q$ queries and 
soundness $s \le \frac{2^{\Omega\(q^{2/3}\)}}{2^q}$.

\noindent {\bf Discussion}

Analyzing acceptance probability of a
low-degree test at distances larger than the unique-decoding radius
seems to require a different set of techniques. It general, to prove
that a code is locally testable, one needs to upper
bound acceptance probability by a function of the distance. This is
achieved by showing that if acceptance probability of the test on an element
$f \in \B^N$ is higher than a certain threshold, there is a
codeword $g$ not far from $f$. In most cases the test itself is used
to efficiently ``decode'' $f$, viewed as a corrupted codeword, to the
unique nearest codeword $g$. This approach is harder to implement when
there are several possible codewords to choose from, and symmetry
breaking in required. The only example we are aware of is the Hadamard
code. In this case one is assisted by 
the fact that the elements of the code are pairwise orthogonal (as
vectors over the reals). In particular, for any $\e > 0$, there could
be only a constant number of codewords at distance smaller than $1/2 -
\e$ from $f$. This no longer holds for degree-$2$ polynomials. For
instance, the list-decoding radius here is $1/4$. Our
main tools in this case are harmonic analysis and additive number
theory. In fact, a significant part of our proof follows the approach
of Gowers \cite{Gowers01} in his proof of Szemeredi's theorem for
arithmetic progressions of length $4$.

\section{Abelian homomorphism testing}
\label{sec:hom}
Let $G$ and $H$ be two finite Abelian groups. In the Abelian
Homomorphism testing problem we are given an oracle 
access to a transformation $\phi~:~G\rarrow H$ and we have to decide
whether $\phi$ is a homomorphism or is at least $\d$-far from any
homomorphism between $G$ and $H$. This problem is a generalization of
the linearity testing problem, in which case $G = H = \Z_2$. It was
first studied in \cite{BLR}, where the following natural
generalization of the basic linearity test was suggested: choose $x,y
\in G$ at random and check whether $\phi(x+y) = \phi(x) +
\phi(y)$. The analysis of this test leads to the following question.

Let $\phi~:~G\rarrow H$
such that the group law for $\phi$ holds with positive probability.
$$
Pr_{x,y\in G} \left(\phi(x) + \phi(y) = \phi(x+y)\right)
\ge \epsilon.
$$
Let $\rho$ be the maximal $\epsilon'$ such that there exists a
homomorphism $\psi$ from $G$ to 
$H$ such that $Pr_x\(\phi(x) = \psi(x)\) \ge \epsilon'$. The question is
whether $\rho$ can be lower bounded in terms of a function of
$\epsilon$ that is independent of $|G|$. In \cite{BLR} this is shown
to be true if $\epsilon > 7/9$. This lower bound on $\epsilon$ is also
necessary \cite{BC}.  

If both $G$ and
$H$ are powers of $\Z_2$, the lower bound on $\epsilon$ was relaxed to   
$\epsilon > 83/128$ \cite{johan+}.

\noindent {\bf Our results}

We show the following theorem to be a simple consequence of two
results \cite{Gowers01, R} in additive number theory.
\thm
\label{thm:hom-test}
Let $p$ be a prime number, and let $\e > 0$.
Let $G$ be a $p$-group of order $r$ and let $H$ be a power of 
$\Z_p$. Let $\phi~:~G\rarrow H$
such that 
$$
Pr_{x,y\in G} \left(\phi(x) + \phi(y) = \phi(x+y)\right)
\ge \epsilon.
$$
Then there exists a homomorphism $\psi~:~G\rarrow H$ such that
$$
Pr_{x\in G} \left\{\phi(x) = \psi(x)\right\} \ge c \cdot r^{-c'} \cdot
\epsilon^{c''} ,
$$
where $c,c',c''$ are absolute constants (independent of the groups
$G,H$).
\ethm
In particular, if both $G$ and $H$ are powers of $\Z_2$, $\rho$ can be 
lower bounded by a function of $\e$, for any $\e > 0$.
In testing terms, this means that the acceptance probability of the
basic test of \cite{BLR} goes to zero as the distance from the code
(the set of all homomorphisms) goes to one.
\section{Tools}
\label{sec:technical}
In this section we discuss the technical tools used in
this paper. We believe these tools, and their connection to recent
results in additive number theory, might be of independent interest.
\subsection{Generalized averages}
Let $H = (V,E)$ be a hypergraph on $t$ vertices. Given a boolean
function $f:~\B^n \rarrow \H$, the 
acceptance probability of the linearity test associated with $H$ on
$f$ is easily seen (cf. Appendix~\ref{app:B}) to be an average of
expressions of the following 
type. Let ${\cal S} = \left\{e_1,\ldots,e_T\right\}$ be a family of edges of
$H$, this is to 
say subsets of $\{1,\ldots,t\}$. We define
the {\it average} of $f$ on ${\cal S}$ in the following way:
\begin{equation}
{\mathbb E}_{{\cal S}}(f) := 
{\mathbb E}_{y_1,\ldots,y_t \in \{0,1\}^n} 
\prod_{j=1}^T f\left(\sum_{i\in e_j}
y_i\right).
\label{matrix_aver}
\end{equation}
The operator $\E_{\cal S}$ is naturally associated with a binary matrix
$A$ whose columns are characteristic vectors of $e_j$. We will also denote this
operator by $\E_A$.
\xam
Let 
$
A = \left[
\begin{array}{ccc}
1 & 0 & 1\\
0 & 1 & 1 
\end{array}
\right]
$.
Then $\E_A(f) = {\E}_{x,y} f(x)f(y)f(x+y)$
is the basic linearity test of \cite{BLR}.
\exam
For $A = [1]$, the average of $f$ over $A$ is, of course, simply the
expectation $\E f$. The notion of generalized average is naturally extended to
real or complex valued functions on $\B^n$.  

The analysis of the probability of acceptance of a hypergraph test
entails studying generalized averages of functions. In particular, we
would like to upper bound such averages by expressions which are
convenient to deal with. 

With this in mind, we define a useful family of binary matrices.

\dfn
\label{dfn:A_k}
For an integer $k\ge 1$ let $A_k$ be a $(k+1)\times 2^k$ matrix of the
following form: the last row of $A_k$ is an all-$1$ vector. Removing
this last row gives a $k\times 2^k$ matrix whose columns are all
binary vectors of length $k$ (in an arbitrary order).
\edfn
Observe that $E_{A_k}(f)$ is precisely $\|f\|^{2^k}_{U_k}$.

We prove several
properties of generalized averages in Appendix~\ref{app:B},
leading to the following main claim. This is essentially a restatement
of theorem~\ref{thm:gowers-norms}.
\thm
\label{thm:test_main}
Assume that all the columns in $A$ are distinct and have at most $k$
ones. Then for any Boolean function $f$  
$$
\Big | \E_A(f) \Big | \le \(\E_{A_k}(f)\)^{\frac{1}{2^k}} = \|f\|_{U_k}
$$
\ethm
\subsection{Gowers norms and pseudorandomness}
In the previous subsection we have seen how Gowers uniformity norms
$\|\cdot\|_{U_k}$ appear naturally in the analysis of linearity
tests. These norms were originally defined in \cite{Gowers01} and
were instrumental in the new proof of Szemeredi's theorem on
arithmetic progressions given in that paper. We refer to
\cite{Gowers01}, \cite{GT05} for a more
detailed discussion. Here let us briefly mention that, intuitively,
the $k$-th uniformity norm of a function is high if this function has
a non-negligible correlation with a polynomial of degree $k-1$. This
is to say, this function has a non-trivial combinatorial structure. On
the other hand, if a function has small uniformity norms, we would
like to deduce that it is 'pseudorandom', in an appropriate sense. In
particular, it is shown in \cite{Gowers01} that if a characteristic
function of a subset of the integers has small uniformity norms, then
the number of arithmetic progressions it contains is similar to that
contained by a random subset of the same size.

This notion of pseudorandomness naturally generalizes (and
strengthens) the standard notion of a boolean function (or a set) being
pseudorandom if its non-zero Fourier coefficients are small. In fact,
the maximal size of a Fourier
coefficient is controlled by the second uniformity norm. Since
uniformity norms are monotone increasing, if a function has a small
$k$-th uniformity norm, $k \ge 2$, it is also pseudorandom in the usual
sense. This is, of course, intuitively clear, since a function far
from degree-$(k-1)$ polynomials is, in particular, far from linear
polynomials. 

In our context, a function $f$ with a small $k$-th uniformity norm, is
pseudorandom in the following sense.
Consider a linearity test associated with a hypergraph $H = (V,E)$ with maximal
edge-size $k$. Theorem~\ref{thm:gowers-norms} implies that the $|E|$ copies
of the basic linearity test that $H$ runs on $f$ behave essentially
independently.  
\subsection{Quadratic Fourier Analysis}
We would now like to give a more specific meaning to
the intuitive notion that a function with a high $k$-th uniformity
norm should have a non-trivial combinatorial structure, presumably a
non-trivial correlation with a polynomial of degree $k-1$.

Unfortunately, at this point, we can only do it for $k=3$. 
By Theorem~\ref{thm:inverse-quadratic}, if $\|f\|_{U_3} \ge \e$
then there exists a quadratic polynomial $g$ such that the distance
between $f$ and $g$ is at most $1/2 - \e'$, for $\e'$ depending on
$\e$ only. 

We conjecture a similar statement to be true for any fixed
$k$. A step in this direction was made in \cite{ST06}, where a function
with a high $k$-th uniformity norm is shown to have variables with
large influence.

Similar results for $k=3$, but replacing $\Z^n_2$ by finite Abelian groups of
cardinality indivisible by $6$, have been independently proved by Green
and Tao \cite{GT05}. The dependence of $\e'$ on $\e$ in both cases is
super-exponential. In \cite{GT05} this dependence is improved in
the following way: it is shown,
specializing here to $\Z_5$ for clarity, that one can find a subspace
$V$ of $\Z^n_5$ of a fixed 
co-dimension and a family of quadratic polynomials $g_y$ indexed by
cosets of $V$, such that typically $f$ is $1/2-\e''$ close to $g_y$ on
$y+V$, where the dependence of $\e''$ on $\e$ is polynomial. This
extension turns out to be useful in obtaining good bounds on
arithmetic progressions of length $4$ in subsets of $\Z^n_5$ (and in
general finite Abelian groups). In this context, Green and Tao introduce the
notion of {\it quadratic Fourier analysis}
\cite{Green-notes}. According to
this point of view, the subject of classical Fourier analysis is
to represent a 
function as a combination of several linear functions (elements of the
Fourier basis) it has non-negligible correlation with (i.e., 
corresponding  Fourier coefficients are large), and of a 'random' remainder
(a function with small Fourier coefficients). In quadratic Fourier
analysis, a function is approximated by a combination of quadratic
polynomials. This approach has proven to be quite effective in
additive number theory \cite{GT:primes,GT05} and in
ergodic theory \cite{Host-Kra, Ziegler}, in situations in which
classical Fourier analysis fails.

Theorems~\ref{thm:inverse-quadratic} and \ref{thm:hypergraph-quadraticity-test}
can be viewed as an application of quadratic Fourier analysis on
$\Z^n_2$ to
boolean functions. We suggest that this tool might have other
applications as well. 
(Among other things, it should be 
possible to extend Theorem~\ref{thm:inverse-quadratic} to obtain
results similar to those of \cite{GT05}, but we haven't checked the details.) 
\newpage

\section{Appendix A: A Proof of Theorem~\ref{thm:inverse-quadratic}} 
\label{app:A}
We start with a short discussion on discrete directional derivatives
of functions on $\B^n$.

Let $f:\B^n \rarrow \H$ be a boolean function, and $y$ be a vector in
$\B^n$. We define the ``derivative of $f$ in direction $y$'' by 
$$
f_y(x) = f(x) f(x+y) 
$$
The transformation $f \mapsto f_y$ is a linear operator. This operator
decreases the 
degree of the polynomial representation of a function: if $f$ is
representable by an $n$-variate 
polynomial of degree $d$, then $f_y$ is representable by a polynomial
of degree $d-1$. 

We define recursively $f_{y_1,y_2} = \(f_{y_1}\)_{y_2}$. It is easy to
see that $f_{y_1,y_2} = f_{y_2,y_1}$, and in fact $f_{y_1,y_2}(x) =
f_{y_2,y_1}(x) = 
f(x) f\(x+y_1\) f\(x+y_2\) f\(x+y_1+y_2\)$. Similarly, the $k$-th
order directional derivative $f_{y_1,...,y_k}$ 
of $f$ with respect to $y_1,...,y_k$ at a point $x$ is given by 
$$
f_{y_1,...,y_k}(x) = \prod_{S\subseteq [k]} f\(x + \sum_{i\in S} y_i\)
$$
If a function $f$ is a polynomial of degree $d$,
then the derivative $f_{y_1,...,y_k}$ is
a polynomial of degree $d-k$, for all choices of linearly independent
$y_1,...,y_k$ \cite{AKKLR, mcw-sloane}. In particular, the
$(k+1)$-th derivative of a degree-$k$ polynomial vanishes (in our terms,
it is identically $1$).

Observe that, in light of the definition above, the claim of
Theorem~\ref{thm:inverse-quadratic} can be interpreted as follows: if
a random third derivative of a function vanishes with probability
greater than $1/2$ then the function is somewhat close to a quadratic
polynomial. 

The proof of the theorem involves several technical lemmas. 
The main tools are Fourier analysis on $\Z^n_2$ (\cite{KKL}) and
additive number theory.
  
In the following the 
Greek letters $\epsilon, \epsilon', \delta, \delta'$ will denote
absolute positive constants 
(independent of $n$) whose value may fluctuate. 
\lem
\label{many_dir}
For a function $f:~\B^n \rarrow \R$ 
$$
\|f\|^8_{U_3} = {\E}_y \sum_{\alpha} \hat{f_y}^4(\alpha)
$$
\elem
\prf
We start with proving
$$
\|f\|^4_{U_2} = \sum_{\alpha} \hat{f}^4(\alpha)
$$
Indeed, 
$$
\|f\|^4_{U_2} = \E_{x,y,z} f(x) f(x+y) f(x+z) f(x+y+z) = \E_x
\(f(x)\cdot \E_{y,z} f(x+y) f(x+z) f(x+y+z)\)
$$
Introducing new variables $u = x+y$, $v = x+z$, this equals to
$$
\E_x \(f(x)\cdot \E_{u,w} f(u) f(w) f(x + u + w)\) = \E_x \(f(x)\cdot
\E_u (f\ast f)(x+u)\) = 
$$
$$
\E_x f(x) (f\ast f\ast f)(x) = \<f,f\ast f\ast
f\> = \<\hat{f}, \hat{f}^3\> = \sum_{\alpha} \hat{f}^4(\alpha)
$$
Now, 
$$
\|f\|^8_{U_3} = {\E}_{x,y,z,w}
f(x)f(x+y)f(x+z)f(x+w)f(x+y+z)f(x+y+w)f(x+z+w)f(x+y+z+w) = 
$$
$$
{\E}_y {\E}_{x,z,w} f_y(x) f_y(x+z) f_y(x+w) f_y(x+z+w) =
\E_y \sum_{\alpha} \hat{f_y}^4(\alpha)
$$
\eprf
\cor
Assuming $f$ is boolean and $\|f\|_{U_3} \ge \epsilon$, 
there exist constants $\delta, \delta'$
and a choice function $\phi:\{0,1\}^n \rarrow \{0,1\}^n$ such
that 
$$
Pr_y\left(\big |\hat{f_y}(\phi(y))\big | \ge \delta\right) \ge \delta'.
$$
\ecor
\prf
The derivatives $f_y$ are also boolean functions, and therefore
$$
\E_y \sum_{\alpha} \hat{f_y}^4(\alpha) \le \E_y \max_{\alpha}
\hat{f_y}^2(\alpha) \cdot \sum_{\alpha} \hat{f_y}^2(\alpha) = \E_y
\max_{\alpha} \hat{f_y}^2(\alpha) 
$$
\eprf

Let $A$ be an $n\times n$ matrix over $\F_2$. If  $g =
(-1)^{\left<Ax,x\right> + a}$ is a quadratic
polynomial\footnote{Observe that, working with the field of $2$
  elements, we can incorporate the linear term of
  a quadratic form in the exponent into the quadratic term, by
  modifying the diagonal of the matrix appropriately.}, then 
$f_y(x) = 
(-1)^{\left<(A + A^t)y,x\right> + a} = (-1)^{\left<By,x\right>
  + a}$. Here $B = A + A^t$ is a symmetric matrix with a zero
diagonal. So for a quadratic polynomial
the choice function $\phi(y) = By$ is linear, and of a special form. 
We will therefore look for similar
properties of the choice function in our case.

It is sufficient to find a choice function that
coincides with an appropriate linear function with positive probability.
This will follow from an observation that if derivatives of two
boolean functions are close on average then so are the functions
themselves (up to a linear shift).
\lem
\label{close_dir_imp_close}
For boolean functions $f,g$: 
$$
{\E}_x \left(\left<f_x,g_x\right>\right)^2 = \sum_{\alpha}
\widehat{fg}^4(\alpha).  
$$
\elem
\prf
$$
{\E}_x \left(\left<f_x,g_x\right>\right)^2 =
{\E}_x {\E}_{y_1,y_2} f_x(y_1) g_x(y_1) f_x(y_2) g_x(y_2) = 
$$
$$
{\E}_x {\E}_{y_1,y_2} f(y_1) f(y_1 + x) g(y_1) g(y_1 + x)
f(y_2) f(y_2 + x) g(y_2) g(y_2 + x) = 
$$
$$
{\E}_x \left({\E}_y (fg)(y) (fg)(y+x) \right)^2 = 
{\E}_x \left((fg)\ast (fg)\right)^2(x) = 
\sum_{\alpha} \widehat{fg}^4(\alpha).  
$$
\eprf

\cor
\label{iso_imp_close_poly}
Let $B$ be a symmetric matrix
with a zero diagonal such that 
$$
{\E}_y \hat{f_y}^2(By) \ge \epsilon.
$$
Then there exists a quadratic polynomial $g$ such that 
$$
\|f - g \| \le \frac12 - \epsilon'.
$$
\ecor
\prf
Let $A$ be a matrix such that $A + A^t = B$. 
Consider a quadratic polynomial $h(x) =  (-1)^{\left<x,Ax\right>}$.
We have
$$
{\E}_x \left(\left<f_x,h_x\right>\right)^2 = {\E}_x
\hat{f_x}^2(Bx) \ge \epsilon'. 
$$
By lemma~\ref{close_dir_imp_close} there is a vector $\alpha$ such that $\big |
\widehat{f h} (\alpha) \big | \ge \epsilon'$. This implies
that there is a choice of $a \in \{0,1\}$ such that for a quadratic
polynomial $g(x) =  (-1)^{\left<x,Ax\right> + \left<x,\alpha\right>
  +a}$ holds 
$$
\|f - g \| \le \frac12 - \epsilon'.
$$
\eprf

We start by finding a weakly linear choice function. This is made
possible by the following observation.
\lem
\label{weak_lin}
$$
{\E}_{x,y} \sum_{\alpha, \beta} \hat{f_x}^2(\alpha)
\hat{f_y}^2(\beta) \widehat{f_{x+y}}^2(\alpha+\beta) = {\E}_y
\sum_{\alpha} \hat{f_y}^6(\alpha). 
$$
\elem
\prf
We start with an observation that for a boolean function $f$ and for
any $x, s$ in $\B^n$ holds $(f_x \ast f_x)(s) = (f_s \ast f_s)(x)$.
Indeed, expanding 
$$
(f_x \ast f_x)(s) = \E_y f_s(y) f_s(x+y) = \E_y f(y) f(y+s) f(x+y)
f(x+y+s).
$$
Define a function $F~:~\B^n \rarrow \{-1,1\}$ by taking $F(y) = f(y)
f(y + (x+s))$. Then the last expression is $(F\ast F)(x) = (F\ast
F)(s)$. Expanding $(f_s \ast f_s)(x)$ we get the same result.

Now,
$$
{\E}_{x,y} \sum_{\alpha, \beta} \hat{f_x}^2(\alpha)
\hat{f_y}^2(\beta) \widehat{f_{x+y}}^2(\alpha+\beta) =
$$
$$
{\E}_{x,y} \sum_{\alpha, \beta} {\E}_{u,u'} f_x(u) f_x(u')
w_{\alpha}(u+u') {\E}_{v,v'} f_y(v) f_y(v') w_{\beta}(v+v') 
{\E}_{z,z'} f_{x+y}(z) f_{x+y}(z') w_{\alpha+\beta}(z+z') =
$$
$$
{\E}_{x,y} {\E}_s {\E}_{u,v,z} f_x(u) f_x(u+s) f_y(v)
f_y(v+s) f_{x+y}(z) f_{x+y}(z+s) =
$$
$$
{\E}_s {\E}_{x,y} \left(f_x\ast f_x\right)(s) \left(f_y\ast
f_y\right)(s) \left(f_{x+y}\ast f_{x+y}\right)(s)  = 
$$
$$
{\E}_s {\E}_{x,y} \left(f_s\ast f_s\right)(x) \left(f_s\ast
f_s\right)(y) \left(f_s\ast f_s\right)(x+y) =
$$
$$
{\E}_s \sum_{\alpha} \hat{f_s}^6(\alpha).
$$
\eprf

\cor
$$
\|f\|_{U_3} \ge \epsilon~~~\Longrightarrow ~~~ {\E}_{x,y}
\sum_{\alpha, 
  \beta} \hat{f_x}^2(\alpha) \hat{f_y}^2(\beta)
\widehat{f_{x+y}}^2(\alpha+\beta) \ge  \epsilon'.
$$
\ecor
\prf
By lemmas~\ref{many_dir} and ~\ref{weak_lin} and Holder's inequality.
\eprf

Define a product distribution on functions 
$\phi:\{0,1\}^n \rarrow \{0,1\}^n$ by taking $Pr(\phi(x) = \alpha) =
\hat{f_x}^2(\alpha)$. The choices for distinct values of $x$ are
independent. Let $\delta = \frac{\epsilon'}{6}$. 
Define a random variable $L$ on this probability space, by taking
$$
L(\phi) = Pr_{x,y} \left\{\phi(x) + \phi(y) = \phi(x+y);~~
\hat{f_x}^2(\phi(x)) \ge \delta;~~...~~\widehat{f_{x+y}}^2(\phi(x+y)) \ge
\delta\right\}. 
$$
\lem
\label{A_exists}
$$
{\E}_{\phi} L(\phi) \ge \frac{\epsilon'}{2}.
$$
\elem
\prf
$$
{\E}_{\phi} L(\phi) = {\E}_{x,y} Pr_{\phi} \left\{\phi(x) +
\phi(y) = \phi(x+y);~~ 
\hat{f_x}^2(\phi(x)) \ge \delta;~~...~~\widehat{f_{x+y}}^2(\phi(x+y)) \ge
\delta\right\} = 
$$
$$
{\E}_{x,y} \sum_{\alpha,\beta~:~\hat{f_x}^2(\alpha) \ge \delta;~~...~~ 
\widehat{f_{x+y}}^2(\alpha+\beta) \ge \delta}
\hat{f_x}^2(\alpha) \hat{f_y}^2(\beta) 
\widehat{f_{x+y}}^2(\alpha+\beta) \ge
$$ 
$$
{\E}_{x,y} \sum_{\alpha,\beta}
\hat{f_x}^2(\alpha) \hat{f_y}^2(\beta) 
\widehat{f_{x+y}}^2(\alpha+\beta) - 3\delta \ge \epsilon' - 3\delta \geq
\frac{\epsilon'}{2}.
$$
\eprf

Take $\phi$ for which $L(\phi) \ge \frac{\epsilon'}{2}$. This is the
choice function we choose. Our goal is to find an appropriate linear
transformation $B:\B^n \rarrow \B^n$ such that $\phi$ and $B$ coincide
on a positive fraction of the domain in which 
$\hat{f_x}^2(\phi(x))\ge \delta$.

We will do this in several steps. In the first step we will find an
affine transformation $x \rarrow Dx + z$ such that ${\E}_x
\hat{f_x}^2(Dx + z) \ge \epsilon'$. Then we will gradually modify this
transformation to obtain a symmetric linear transformation $B$
with a zero diagonal such that ${\E}_x
\hat{f_x}^2(Bx) \ge \epsilon'$. By lemma~\ref{iso_imp_close_poly} this
will conclude the proof of the theorem.

The first step is the hardest. We will follow an approach of Gowers
from his proof of Szemeredi's theorem for arithmetic progressions of length
four \cite{Gowers01}. Note that a structural theorem of Freiman for sets
with small sumsets in $\Z$ is replaced by a theorem of Ruzsa for such
sets in $\Z^n_2$.
 
Let $A = \left\{x~:~\hat{f_x}^2(\phi(x)) \ge \delta\right\}$. Then, by
the choice of $\phi$, the cardinality $m$ of $A$ is a positive
fraction of $2^n$, 
and there are $\Omega\left(m^2\right)$ triples $(x,y,x+y)$ in $A^3$
satisfying $\phi(x) + \phi(y) = \phi(x+y)$. 

Now define a subset ${\cal A}$ of $\{0,1\}^{2n}$ as
$$
{\cal A} = \left\{(x,\phi(x))~:~x\in A\right\}.
$$
${\cal A}$ is the graph of $\phi$ on $A$. We have $|{\cal A}| = |A| =
m$, and there are $\Omega\left(m^2\right)$ triples $(a,b,a+b)$ in
${\cal A}^3$. 

\thm (Gowers \cite{Gowers01})
\label{Gowers}
For any subset ${\cal A}$ of an abelian group satisfying above, there
a subset ${\cal A}'$ of ${\cal A}$ containing a constant fraction of
the elements, such that 
$$
|{\cal A}' + {\cal A}'| \le c \cdot |{\cal A}'|,
$$
for an absolute constant $c$.
\ethm
\thm (Ruzsa \cite{R})
\label{Ruzsa}
Let $G$ be an abelian group. Assume that the order of the elements in
$G$ is bounded, and let $r$ be the maximal order of an element.
Let ${\cal A}' \subseteq G$ with the property above. Then
$$
\frac{|{\cal A}'|}{|<{\cal A}'>|} \ge c'\cdot r^{-c''},
$$
for some absolute constants $c',c''$.
\ethm
We assume the projection of $<{\cal A}'>$ on the first $n$
coordinates to be of full rank. (Otherwise add a finite number of
vectors to $<{\cal A}'>$ to ensure this.) Therefore, there are $n$
vectors $v_1...v_n$ in $\{0,1\}^n$ such that the vectors $u_i =
(e_i,v_i)$ are in  $<{\cal A}'>$. Let $U = <u_1...u_n>$. Clearly $U =
\left\{(x,Dx)~:~x\in \{0,1\}^n\right\}$, where the matrix $D$ is
defined by $De_i = v_i$, $i=1...n$.  
$U$ is a subspace of $<{\cal A}'>$ of a finite co-dimension. Therefore
there exists a vector $c \in \{0,1\}^{2n}$ such that a constant
fraction of the vectors in ${\cal A}$ sit in $U + c$. This is the same
as to say that there is a vector $z \in \{0,1\}^n$ such that for
$\Omega\left(2^n\right)$ points $x \in {\cal A}$ holds $\phi(x) = Dx +
z$.
Alternatively:
$$
{\E}_x \hat{f_x}^2(Dx + z) \ge \epsilon'.
$$
We can choose $z=0$.
\lem
\label{positive_dfn}
Define a function $F~:~\{0,1\}^n
\rarrow {\R}$ by $F(z) = \sum_y \hat{f_y}^2(Dy + z)$. Then
$$
\hat{F}(x) = \hat{f_x}^2(D^t x).
$$
\elem
\prf
$$
\hat{F}(x) = {\E}_z F(z) w_x(z) = {\E}_z w_x(z) \sum_y
\hat{f_y}^2(Dy + z) = 
$$
$$
{\E}_z w_x(z) \sum_y {\E}_{u,u'} f_y(u) f_y(u') w_{Dy+z}(u+u') = 
{\E}_{y,u} f_y(u) f_y(u+x) w_{Dy}(x) =
$$
$$
{\E}_y (f_y\ast f_y)(x)  w_{Dy}(x) = {\E}_y (f_x\ast f_x)(y)
w_{D^t x}(y) = \hat{f_x}^2(D^t x).
$$
\eprf

Since the transform of $F$ is nonnegative, $F$ attains its maximum in
$0$. Therefore 
$$
{\E}_x \hat{f_x}^2(Dx) \ge {\E}_x \hat{f_x}^2(Dx + z) \ge \epsilon'.
$$
We want to replace $D$ by a symmetric matrix. The following fact is useful.
\lem
\label{zero_on_not_ort}
$$
\hat{f_y}(x) = 0
$$
for any $x$ and $y$ with $\left<x,y\right> = 1$.
\elem
\prf
$$
\hat{f_y}(x) = {\E}_z f_y(z) w_x(z) =  {\E}_z f(y) f(z) f(y+z)
w_x(z) =  {\E}_z f(y) f(z) f(y+z) w_x(y+z) = 
$$
$$
w_x(y) {\E}_z
f(y) f(z) f(y+z) w_x(z) =  w_x(y) \hat{f_y}(x) = - \hat{f_y}(x) 
$$
\eprf

Therefore, for $g(x) = (-1)^{\left<x,Dx\right>}$ holds
$$
{\E}_x g(x) \hat{f_x}^2(Dx) = {\E}_x \hat{f_x}^2(Dx) \ge
\epsilon'.
$$
On the other hand
$$
{\E}_x g(x) \hat{f_x}^2(Dx) = \sum_z \hat{g}(z) {\E}_x
\hat{f_x}^2(D^t x + z).
$$
Since the numbers $\lambda_z = {\E}_x \hat{f_x}^2(D^t x + z)$ are
nonnegative and sum to one, we deduce by Jensen's inequality that 
$$
\sum_z \hat{g}^2(z) {\E}_x \hat{f_x}^2(D^t x + z) \ge
(\epsilon')^2.
$$
However
$$
\sum_z \hat{g}^2(z) {\E}_x \hat{f_x}^2(D^t x + z) = 
{\E}_x (g\ast g)(x) \hat{f_x}^2(Dx) = 
{\E}_x \delta_{Dx,D^t x}\cdot g(x) \hat{f_x}^2(Dx).
$$
Let $S$ be a matrix defined in the following way: Set $U =
\left\{x~:~Dx =  D^t x\right\}$. Then $U$ is a subspace of $\Z^n_2$. Let $S$ be
defined on $U$ by taking $S(x) = D(x)$ on $U$. Then for any $x,y\in U$
holds $\left<x, Sy\right> = \left<Sx, y\right>$. Now the definition of
$S$ could be extended to the whole space keeping this
property. Therefore $S$ is a symmetric matrix such that
$$
{\E}_x  \hat{f_x}^2(Sx) \ge \epsilon.
$$
It remains to deal with the diagonal of $S$. Let $f$ be the vector on
the diagonal of $S$. Since $\left<x,Sx\right> = \left<x,f\right>$, we have 
$$
\frac{1}{2^n} \sum_{x \perp f} \hat{f_x}^2(Sx) \ge \epsilon.
$$
Define a matrix $B$ by taking $Bx = Sx$ if $x \perp f$, and extending
$B$ appropriately to the whole space. Namely for $w \not \perp f$
take $Bw = z$, so that
$\left<z,x\right> = \left<w,Bx\right> = \left<w,Sx\right>$ for all $x
\perp f$, and $\left<w,z\right> = 0$.
Then $B$ is symmetric with zero diagonal, and 
$$
{\E}_x  \hat{f_x}^2(Bx) \ge \epsilon'.
$$
This concludes the proof of theorem~\ref{thm:inverse-quadratic}, but
for the dependence of $\e'$ on $\e$. Tracing this dependence through
the proof, it is possible to see that we
can choose $\e' \ge
\Omega\(\mbox{exp}\left\{-\(\frac{1}{\e^C}\)\right\}\)$ 
for an absolute constant $C$. 
\eprf

\section{Appendix B: A Proof of Theorem~\ref{thm:gowers-norms}}
\label{app:B}
We will prove Theorem~\ref{thm:test_main}. 

This will imply
Theorem~\ref{thm:gowers-norms} as follows: Let $f:~\B^n \rarrow \H$ be
a boolean function. Let $H = (V,E)$ be a hypergraph with maximal
edge-size $d$. For a subset ${\cal S} \subseteq E$ of edges of $H$,
let $\sigma(S) = \sum_{e\in {\cal S}} (|e|+1)$. 
The probability that the linearity test associated with
$H$ accepts $f$ is given by 
$$
\frac{1}{2^{|E|}}
\sum_{{\cal S} \subseteq E} f^{\sigma({\cal S})}(0) \cdot
\E_{x_1,\ldots,x_t}
\left[ \prod_{e \in {\cal S}} \prod_{i\in e} f(x_i) \cdot
f\left(\sum_{i\in e} x_i\right) \right] 
$$
The summand corresponding to $S = \emptyset$ is $1$. By
Theorem~\ref{thm:test_main} all the other
summands are at most $\|f\|_{U_d}$ in absolute value. 
Theorem~\ref{thm:gowers-norms} follows.

We start with some facts on generalized averages
(\ref{matrix_aver}). Recall that each 
such average is naturally associated with a $t\times T$ binary matrix $A$. 
The first observation is that some matrices (families of sets)
define the same average operator.
\lem
\label{holistic} 
Multiplying $A$ on the left by a non-singular
$t\times t$ matrix $B$ does not change the value of the average, namely
for any function $f$ holds ${\mathbb E}_{A}(f) = {\mathbb E}_{BA}(f)$,
\elem
\prf
$\E_A$ and $\E_{BA}$ are the same up to order of summation.  
\eprf
\cor
We may (and will) assume that the rows of $A$ are linearly
independent, since if $rank(A) = r < t$ we can choose a 
non-singular $t\times t$ matrix $B$, so that in $BA$ the last $t-r$ rows are
zeroes, and consequently can be removed without changing the value of 
$\E_A$. 
\ecor

Consider an equivalence relation on $t\times T$ binary matrices of
full rank,
defined by left multiplication by a non-singular $t\times t$
matrix. It is easy to see that two matrices are equivalent iff 
their rows span the same $t$-dimensional space over $\Z^T_2$
(or they represent the same rank-$t$ {\it binary matroid} on
$\{1...T\}$ \cite{ox}). 

The following definition and lemmas are natural (and well-known)
in the setting of matroids (see \cite{ox}). 
\dfn
A {\it hyperplane} of $A$ is a maximal subset of $\{1,\ldots,T\}$
such that the columns of $A$ indexed by this subset are {\it not} of
full rank. 
\edfn
\lem
\label{all_same}
A vector $v$ is a minimal non-zero vector in the row space of $A$ iff
the complement of its support is a hyperplane of $A$.
\elem

\lem
The row space of $A$ is spanned by its minimal  non-zero vectors.
\elem

The key part of the proof of theorem~\ref{thm:test_main} 
is the following technical proposition:
\pro
\label{step}
Let $A$ be a full rank $t\times T$ binary matrix. 
Let $v$ be a minimal vector in the row-space of $A$. Let $A'$ be 
a $(t+1)\times 2|v|$ matrix obtained
from $A$ by the following procedure:
\begin{enumerate}
\item
Delete all the columns not in the support of $v$, obtaining a $t\times
|v|$ matrix $B$. 
\item 
Set
$$
A' = 
\left[
\begin{array}{c|c}
\hspace*{.3in} B \hspace*{.3in} & \hspace*{.3in} B \hspace*{.3in} \\
1\ldots 1                       & 0\ldots 0
\end{array}
\right]
$$
\end{enumerate}
Then for all boolean functions $f$, 
$$
{\mathbb E}_{A}(f) \le \sqrt{{\mathbb E}_{A'}(f)}.
$$
\epro

\xam
Let
$$
A = \left[
\begin{array}{ccc}
1 & 0 & 1\\
1 & 1 & 0
\end{array}
\right],
$$
and take $v = \{1,3\}$. Then 
$$
A' = \left[
\begin{array}{cccc}
1 & 1 & 1 & 1\\
1 & 0 & 1 & 0\\
1 & 1 & 0 & 0
\end{array}
\right].
$$
\exam

\prf
Let $H$ be the complement of $v$, and by lemma~\ref{all_same} a
hyperplane of $A$. 
Let $H_0$ be a maximal independent subset (a basis) of
$H$, of size $t-1$. We assume that $H = \{1,\ldots, |H|\}$ and 
$H_0 = \{1,\ldots,t-1\}$. Multiply $A$ on the left by a non-singular
$t\times t$ matrix $B$ so that the first $t-1$
columns of $A$ are the first $t-1$ unit vectors. Since $H_0$ is a
basis of $H$, the columns of $BA$ indexed by $H$ will have a zero in
their last coordinate, while the columns in $v = H^c$ will have one (since
$H$ is a hyperplane). Namely
$$
BA = 
\left[
\begin{array}{c|c}
10\ldots 0 & \hspace{.5in}\\
01\ldots 0 & \hspace{.5in}\\
\vdots     & \mbox{\Huge N}\\
00\ldots 1 & \hspace{.5in}\\
\hline
00\ldots 0 & 0\ldots 0 1\ldots 1
\end{array}
\right].
$$
We have, for any $f$:
$$
{\mathbb E}_A(f) = {\mathbb E}_{BA}(f) = {\mathbb E}_{y_1,\ldots,y_t} 
f(y_1)\cdot\ldots \cdot f(y_{t-1}) \cdot \prod_{i=t}^T
f\left(\prod_{j\in A_i} y_j\right).
$$
We upper bound the right hand side in the following way, applying the
Cauchy-Schwarz inequality:
$$
{\mathbb E}_M(f) = {\mathbb E}_{y_1,\ldots,y_{t-1}} 
F\left(y_1,\ldots,y_{t-1}\right) \cdot
G\left(y_1,\ldots,y_{t-1}\right) \le
$$
$$
\sqrt{{\mathbb E}_{y_1,\ldots,y_{t-1}}
F^2\left(y_1,\ldots,y_{t-1}\right)} \cdot
\sqrt{{\mathbb E}_{y_1,\ldots,y_{t-1}}
G^2\left(y_1,\ldots,y_{t-1}\right) },
$$
where $F\left(y_1,\ldots,y_{t-1}\right) = 
\prod_{j=1}^{t-1} f(y_j)$, and $G\left(y_1,\ldots,y_{t-1}\right) = 
{\mathbb E}_{y_t} \prod_{i=t}^T f\left(\prod_{j\in A_i} y_j\right)$.

Observe that
$$
{\mathbb E}_{y_1,\ldots,y_{t-1}} F^2\left(y_1,\ldots,y_{t-1}\right) = 
{\mathbb E}_{y_1,\ldots,y_{t-1}} f^2(y_1)\cdot\ldots\cdot f^2(y_{t-1})
= {\mathbb E} 1 = 1.
$$
Therefore ${\mathbb E}_M(f) \le \sqrt{{\mathbb E} G^2}$.
It is easily seen that ${\mathbb E} G^2$ presents an average of $f$
on a $(t+1)\times 2(T-t+1)$ matrix $A^{(1)}$, where
$$
A^{(1)} = 
\left[
\begin{array}{c|c}
\hspace{.5in} & \hspace{.5in}\\
\mbox{\Huge N}& \mbox{\Huge N}\\
\hline
0\ldots 0 0\ldots 0 & 0\ldots 0 1\ldots 1\\
\hline
0\ldots 0 1\ldots 1 & 0\ldots 0 0\ldots 0
\end{array}
\right].
$$
We now transform the matrix $A^{(1)}$ to $A'$ in three steps. 
Let $u$ and $v$ be the last two
rows of the matrix. 
First, replace $u$ with $u+v$ obtaining a new matrix
$A^{(2)}$. 

The second step uses booleanity of $f$.
Note that for any matrix $A$ and any boolean function
$f$, deleting a pair of identical columns of $A$ does not change the
average of $f$ on $A$. We delete all the columns of $A^{(2)}$ in which
the last two coordinates are zero, obtaining a $(t+1)\times 2|v|$
matrix $A^{(3)}$. The third step is to multiply $A^{(3)}$ by a
$(t+1)\times (t+1)$ matrix 
$
S_1 = \left[
\begin{array}{cc}
B^{-1} & 0 \\
0 & 1 
\end{array}
\right]
$,
obtaining $A'$.
\eprf

Now we are ready to prove Theorem~\ref{thm:test_main}.
We will prove the theorem for $k = 3$. The proof for larger values of
$k$ is similar.

Let $A$ be a matrix with at most $3$ ones in each column. Assume the
rows of $A$ to be independent. Let $u$ be the first row. The first
step is to replace $u$ with a minimal non-zero vector $v$ with a smaller
support. If $u$ is already a minimal vector let $v = u$. Otherwise there
is a vector $w$ in a row-space of $A$ whose support is strictly
smaller than that of $u$. One of the vectors 
$w$ or $u+w$ is not spanned by the rest of the rows of $A$ and we set
$u_1$ to be this vector. If $u_1$ is minimal set $v =
u_1$. Otherwise continue with $u_1$ instead of $u$. Clearly this
process stops after a finite number of steps and does not change the
row space of $A$.

Now we apply the transformation of proposition~\ref{step} 
to the new matrix $A$, choosing $v$ as the appropriate
minimal vector. Consider the submatrix $B$ of the new matrix $A'$. The
first row of $B$, and therefore of $A'$ as well, is a
$1$-vector. Moving it to be the last row, we obtain 
$$ 
A' = 
\left[
\begin{array}{c|c}
\hspace*{.3in} B' \hspace*{.3in} & \hspace*{.3in} B' \hspace*{.3in} \\
1\ldots 1                       & 0\ldots 0 \\
1\ldots 1                       & 1\ldots 1
\end{array}
\right]
$$
The matrix $B'$ has at most $2$ ones in each column. Note that all the
columns in $B'$ are distinct, since so were the columns of $A$. There
are two cases to distinguish. 

$B'$ has only one column. Then, removing dependent rows, we get to a
$2\times 2$ matrix 
$
A_1 = \left[
\begin{array}{cc}
1 & 0 \\
1 & 1 
\end{array}
\right]
$. 
By Proposition~\ref{step}, for any boolean function $f$ holds
$
\Big | \E_A(f) \Big | \le \(\E_{A_1}(f)\)^{1/2} = \|f\|_{U_1} \le \|f\|_{U_3}.
$
In the last inequality we use monotonicity of uniformity norms.

$B'$ has more than one column. If there are dependencies between rows
of $A'$ we remove them, keeping only a spanning set of rows, starting
with the last two. In particular $A'$ has no all-$1$ rows except the last.

We now repeat
the procedure starting from $A'$. Consider the first row $u$ of
$A'$. $u = (u'~|~u')$ is a symmetric vector. 
If it is minimal set $v = u$. If not, there is a minimal vector $w$
of smaller support such that replacing $u$ with $w$ does not affect
the row space of $A'$. The vector $w$ is in this row space, therefore it is
either symmetric or antisymmetric. However if it is antisymmetric then
$u$ has to be an all-$1$ row, which we have excluded. Therefore we can
replace $u$ by a symmetric minimal vector $v$. Now apply the proposition
with $A'$ and $v$, and obtain a new matrix (after simplification)
$$
A'' = 
\left[
\begin{array}{c|c|c|c}
\hspace*{.3in} B'' \hspace*{.3in} & \hspace*{.3in} B'' \hspace*{.3in}
& \hspace*{.3in} B'' \hspace*{.3in} & \hspace*{.3in} B'' \hspace*{.3in} \\
1\ldots 1   & 0\ldots 0 & 1\ldots 1   & 0\ldots 0 \\
1\ldots 1   & 1\ldots 1 & 0\ldots 0   & 0\ldots 0 \\
1\ldots 1   & 1\ldots 1 & 1\ldots 1   & 1\ldots 1
\end{array}
\right]
$$
The matrix $B''$ has at most one $1$ in each column. All the
columns in $B'$ are distinct.

Once again, there are two cases. If there is only one column in $B''$,
after simplification we get a $3\times 4$ matrix
$$
A_2 = \left[
\begin{array}{cccc}
1 & 0 & 1 & 0 \\
1 & 1 & 0 & 0 \\
1 & 1 & 1 & 1 
\end{array}
\right]
$$
such that for any boolean function $f$ holds $\Big |\E_A(f) \Big | \le
\(\E_{A_2}(f)\)^{1/4} = \|f\|_{U_2} \le \|f\|_{U_3}$. 

If $B''$ has more than one column we iterate once again. It is not
hard to see that the new matrix $B'''$ will necessarily have only one
column. After simplifying, we will get to
a $4\times 8$ matrix
$$
A_3 = \left[
\begin{array}{cccccccc}
1 & 0 & 1 & 0 & 1 & 0 & 1 & 0\\
1 & 1 & 0 & 0 & 1 & 1 & 0 & 0\\
1 & 1 & 1 & 1 & 0 & 0 & 0 & 0\\
1 & 1 & 1 & 1 & 1 & 1 & 1 & 1 
\end{array}
\right]
$$
such that for any boolean function $f$ holds $\Big | \E_A (f) \Big |
\le \(\E_{A_3}(f)\)^{1/8} = \|f\|_{U_3}$. The theorem is proved. 
\eprf
\section{Appendix C: Other proofs}
\label{app:C}
\subsection{Proof of Proposition~\ref{pro:dicho}}
First we estimate $\|f\|_{U_3}$ within additive precision of
$O(\d)$. This can be done by choosing at random
$\Omega\(\frac{1}{\d}\)$ quadruples of vectors $x,y_1,y_2,y_3 \in
\B^n$ and averaging $f_{y_1,y_2,y_3}(x)$ over the choices. Let us call
this average $\nu$. It is easy to see that for a sufficiently large
number of sampled quadruples, $\Big | \nu - \|f\|_{U_3} \Big | \le \d/2$ with
high probability.   

Assuming this is true, there are two possibilities. First, $\nu \ge
\d$. In this case, $\|f\|_{U_3} \ge
\d/2$. Theorem~\ref{thm:inverse-quadratic} now implies that the second
option of the proposition holds. 

The second option is $\nu < \d$. In this case $\|f\|_{U_3} <
3\d/2$. Now we follow an argument from \cite{GT05}. Let $g$ be a 
quadratic polynomial. Recalling the interpretation of $\|\cdot\|^8_{U_3}$
as the average of the third order derivative of a function, it is easy
to see that $\|fg\|_{U_3} = \|f\|_{U_3} \le 3\d/2$. Observe that
the first uniformity ``norm'' of a function is the square of its
expectation. By the monotonicity of uniformity norms, 
$$
\<f,g\> = \widehat{fg}(0) \le \|f\|^{1/2}_{U_1} \le \|f\|^{1/2}_{U_3}
\le O\(\d^{1/2}\)
$$
Since both $f$ and $g$ are boolean functions, this implies that the
distance between $f$ and $g$ is at least $1/2 - \Omega\(\d^{1/2}\)$,
and the first option of the proposition holds.
\subsection{Proof of Theorem~\ref{thm:hypergraph-quadraticity-test}}
The completeness of the test follows from the fact that it checks whether
third order derivatives of the function vanish.

The fact that the soundness of the test is $1/2^{|E|}$ is an immediate
consequence of Theorems~\ref{thm:test_main} with $k =
3$ together with Theorem~\ref{thm:inverse-quadratic}. Indeed, let a
boolean function $f$ be $1/2 - \e$ far from quadratic polynomials.
Similarly to the proof of Theorem~\ref{thm:gowers-norms}, the
acceptance probability of the test on a function $f$ is upper bounded
by $1/2^{|E|} + \|f\|_{U_3}$. By Theorem~\ref{thm:inverse-quadratic},
this is at most $1/2^{|E|} + \e'$, with $\e' \rarrow 0$ with $\e$.
\subsection{Proof of Theorem~\ref{thm:hom-test}} 
Combining Theorems~\ref{Gowers} and \ref{Ruzsa} similarly to the proof
of Theorem~\ref{thm:inverse-quadratic}, we obtain the following
claim.
\thm
Let $p$ be a prime number, and let $\e > 0$.
Let $G$ be a $p$-group of order $r$ and let $H$ be a power of 
$\Z_p$. Let $\phi~:~G\rarrow H$
such that 
$$
Pr_{x,y\in G} \left(\phi(x) + \phi(y) = \phi(x+y)\right)
\ge \epsilon.
$$
Then there exists a homomorphism $\psi~:~G\rarrow H$ and an element $h
\in H$ such that
$$
Pr_{x\in G} \left(\phi(x) = \psi(x) + h\right) \ge c \cdot r^{-c'} \cdot
\epsilon^{c''} ,
$$
where $c,c',c''$ are absolute constants (independent of the groups
$G,H$).
\ethm
The following lemma concludes the proof of Theorem~\ref{thm:hom-test}. 
\lem
Let $G$ be a $p$-group of order $r$, and let $H$ be a power of 
$\Z_p$. Let $\phi~:~G\rarrow H$ be such that
there exists a homomorphism $\psi~:~G\rarrow H$ and an element $h
\in H$ such that
$$
Pr_{x\in G} \left(\phi(x) = \psi(x) + h\right) \ge \delta.
$$
Then there exists a homomorphism $\psi'~:~G\rarrow H$ such that
$$
Pr_{x\in G} \left(\phi(x) = \psi'(x)\right) \ge \frac{c}{r} \cdot
\delta,
$$
for an absolute constant $c$.
\elem
\prf
Let $E = \{x\in G:~\phi(x) = \psi(x) + h\}$. Let $G = \prod_{i=1}^m
\Z_{p^{k_i}}$. There 
exists an absolute constant $c$, 
a coordinate $1 \le i \le m$, and a generating element $g \in
\Z_{p^{k_i}}$ such  
that for at least $\frac{c}{r}$-fraction of the elements of $E$ holds
$x_i = g$. Call this set $E'$.
Let $e_1...e_m$ be the standard basis of $G$. Consider a
homomorphism $\psi'~:~G\rarrow H$ defined as follows: $\psi'(e_j) =
\psi(e_j)$ for $j \not = i$ and $\psi'(g\cdot e_i) = \psi(e_i) + h$. Then
$\psi'$ agrees with $\phi$ on $E'$.
\eprf

\section{Acknowledgements}
Part of this work was done while visiting Johan H{\aa}stad at KTH. I am
very grateful to Johan for his hospitality and very helpful
discussions. I would also like to thank Michael Ben-Or, Nati Linial,
Luca Trevisan, and Benjamin Weiss for valuable conversations.

\end{document}